\documentclass[11pt,a4paper]{article}
\usepackage[margin=2cm, bottom=3cm]{geometry}
\RequirePackage[round, authoryear]{natbib}
\usepackage[labelfont=bf]{caption}
\usepackage[utf8]{inputenc}
\usepackage[T1]{fontenc}
\usepackage{array}
\usepackage{multirow}
\usepackage{amsmath}
\usepackage{amssymb}
\usepackage{stmaryrd}
\usepackage{graphicx}
\usepackage{wasysym}
\usepackage{diagbox}
\usepackage{amsthm}
\usepackage{enumerate}
\usepackage[shortlabels]{enumitem}
\usepackage[colorlinks,
linkcolor=blue,
anchorcolor=blue,
citecolor=blue]{hyperref}
\usepackage{float}
\usepackage{graphicx}
\usepackage{placeins}
\usepackage{makecell}
\usepackage{xr}
\usepackage{cases}
\usepackage{mathtools}
\usepackage[math]{cellspace}
\usepackage{bigints}
\usepackage{setspace}
\usepackage{enumitem}
\usepackage{times, cleveref}

\theoremstyle{definition}

\title{ 
	\textbf{A new method for generalizing non-self-intersecting flexible polyhedra}
}
	
\author{Zeyuan He, Simon D. Guest 
	\\ \small Department of Engineering, University of Cambridge
	\\ \small zh299@cam.ac.uk, sdg@eng.cam.ac.uk}
	
\date{}

\begin{document}
	
\maketitle

\begin{abstract}
A surface is considered flexible if it allows a continuous deformation that preserves both metric and smoothness. We introduce a novel construction method, called `base + crinkle,' for generating a broad class of non-self-intersecting flexible closed polyhedral surfaces (i.e. flexible polyhedra). These flexible polyhedra can be non-triangulated, exhibit multiple kinematic degrees of freedom, and possess topologies beyond the sphere. The geometric result provides fresh insights into the geometry of origami and the design of engineering mechanisms, such as sealed-chamber robotics and distortion-free metamorphic grippers.
\end{abstract}


\section{Introduction}
\label{sec:introduction}

\begin{figure} [tbph]
	\noindent \begin{centering}
		\includegraphics[width=1\linewidth]{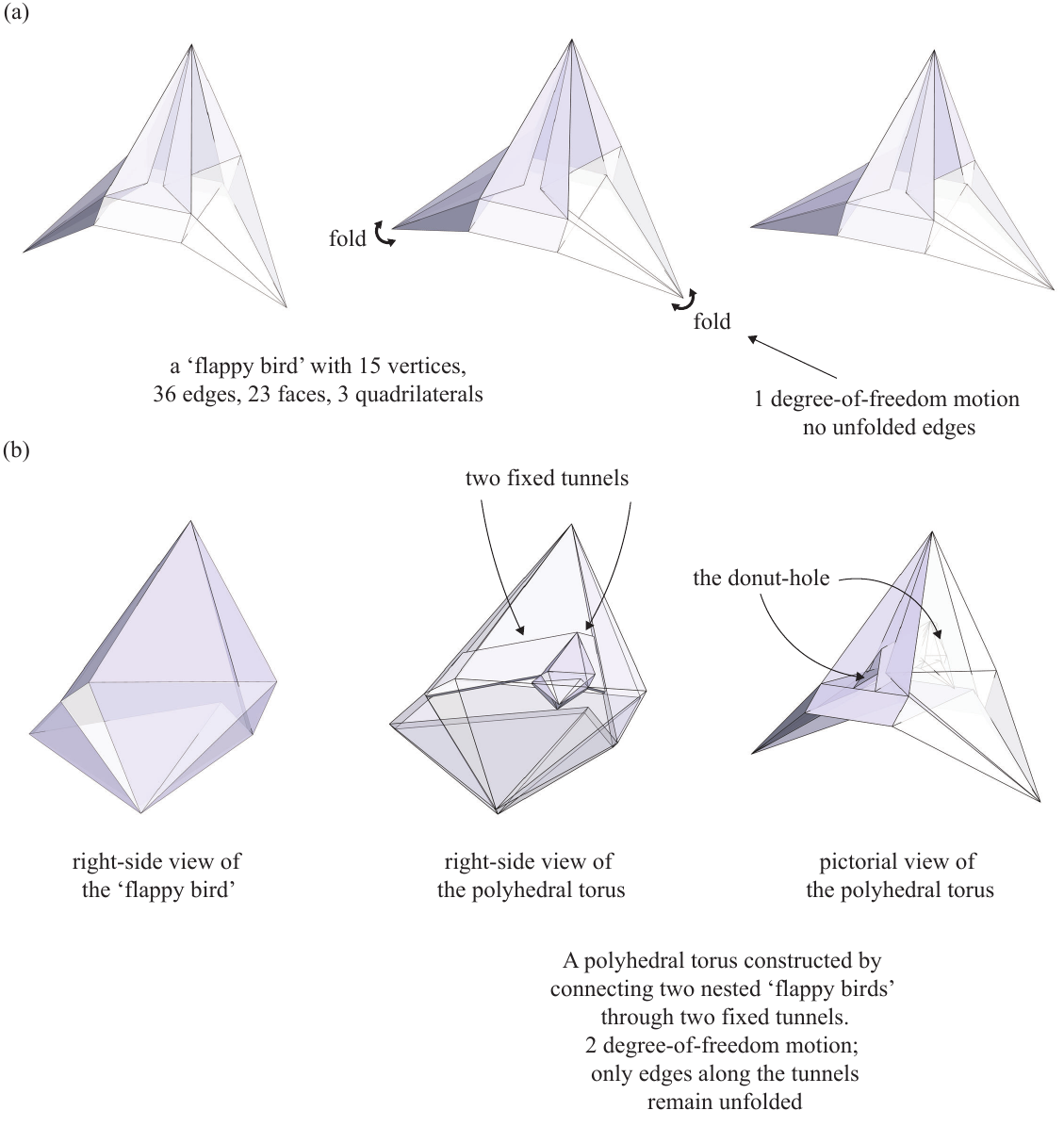}
		\par \end{centering}
	
	\caption{\label{fig: bird}Two new flexible polyhedra generated using the `base + crinkle' method. (a) A non-triangulated, topologically spherical flexible polyhedron exhibiting a single kinematic degree of freedom. It achieves a wider folding range (29.2$^\circ$) compared to the Steffen’s polyhedron (27$^\circ$, Figure \ref{fig: Steffen}) and, notably, contains no unfolded edges --- all dihedral angles vary continuously during its motion. (b) A non-triangulated, topologically toroidal flexible polyhedron with two kinematic degrees of freedom. The construction begins with two nested `flappy birds,' whose folding motions are independent and non-interfering under appropriate parameter adjustments, ensuring no self-intersection occurs. To achieve toroidal topology, two triangles from the inner `flappy bird' are removed, while the fixed tunnels and the triangulated openings on the outer surface together form the handle. Further details of the construction are provided in Section \ref{section: result and method}.}
\end{figure}

Non-self-intersecting (embedded) flexible polyhedra—closed polyhedral surfaces that deform continuously while keeping all faces rigid—are classically known only through a few isolated examples with highly restricted geometries: they are typically triangulated, topologically spherical, and possess a single kinematic degree of freedom (e.g. the Steffen's polyhedron, Figure \ref{fig: Steffen}). Here we introduce a `base + crinkle' construction that substantially broadens this landscape, producing non-triangulated, multi-degree-of-freedom, and topologically toroidal flexible polyhedra (Figure \ref{fig: bird})—and even examples with no unfolded edges—while preserving non-self-intersection. Although not exhaustive, this study demonstrates the potential and versatility of the proposed method through several new examples exhibiting large range of motions, highlighting new avenues for engineering design inspired by these generalized flexible polyhedra.

In this article, a \textit{polyhedron} (i.e. closed polyhedral surface) is a piecewise-linear surface in $\mathbb{R}^3$ obtained by gluing finitely many planar polygons along their edges so that the result is a compact two-manifold without boundary. Such a polyhedron is \textit{self-intersected} (immersed) if its faces overlap or penetrate one another; otherwise, it is non-self-intersecting (embedded). Unless otherwise stated, all polyhedra mentioned in this article are non-self-intersecting. A surface is \textit{flexible} if it admits a continuous deformation that preserves both its intrinsic metric (distance defined on the surface) and smoothness. Polyhedra can have different topologies, such as spherical (genus 0), toroidal (genus 1), double-toroidal (genus 2), and so on. 

Since the time of \citet{euler_opera_1776} and \citet{maxwell_transformation_1819}, researchers have been captivated by the question of whether non-self-intersecting
closed surfaces can exhibit flexibility. Cauchy’s celebrated rigidity theorem \citep{cauchy_sur_1813} provided an early and influential answer—showing that all convex polyhedra are necessarily rigid. Yet, the possibility of flexibility re-emerged when \citet{bricard_memoire_1897} discovered three classes of self-intersecting, non-convex flexible octahedra, highlighting that the challenge was not just flexibility but also the elimination of self-intersections. Remarkably, \citet{connelly_counterexample_1977} constructed the first example of a non-self-intersecting flexible polyhedron using special geometric insights and so disproved the rigidity conjecture of Euler after two-hundred years. Subsequently, \citet[9 vertices]{steffen_symmetric_1978} (Figure \ref{fig: Steffen}) and Kuiper and Pierre \citep[11 vertices]{connelly_rigidity_1979} presented a few other variations, though these remained isolated examples built on similar construction principles. Steffen’s polyhedron contains one unfolded edge along which the dihedral angle remains constant during its motion, which is necessary in his construction to eliminate self-intersection. In her thesis \citet{lijingjiao_flexible_2018}, she introduced an $n$-tetrahedra ($n \ge 2$) variation based on Steffen’s polyhedron, featuring $6n-2$ vertices, $n-1$ kinematic degrees of freedom, and 2 unfolded edges. More recently, \citet{gallet_pentagonal_2024} discovered an 8-vertex flexible polyhedron—the smallest possible—featuring three unfolded edges arranged as a `tetrahedral tent' to prevent self-intersection. Subsequent developments include the 26-vertex construction by \citet{alexandrov_embedded_2024} with no unfolded edges, and two further designs introduced by \citet{atlason_constructing_2025, atlason_cutting_2025}: the 15-vertex `fox,' derived from their twinning method, and a new 8-vertex model offering a substantially larger folding range than that of \citet{gallet_pentagonal_2024}.  These reported examples are all triangulated and topologically spherical. No examples have been documented with non-triangular faces or topologies beyond the sphere. 

From an engineering perspective, a flexible polyhedron represents a largely unexplored class of mechanisms, notable for its intricate geometry and unconventional kinematic behaviour. Certain geometric quantities, such as the enclosed volume \citep{sabitov_problem_1995}, remain invariant throughout its motion. These distinctive properties motivate the generalization of such surfaces and their exploration for geometric adaptability and form-finding applications, including sealed-chamber robotics and distortion-free metamorphic grippers.

\section{Result and method} \label{section: result and method}

We present a novel construction approach, termed the `base + crinkle' method, for producing a broad family of non-self-intersecting flexible polyhedra. The resulting structures may be non-triangulated, support multiple kinematic degrees of freedom, and exhibit topologies beyond that of the sphere. An initial result is shown in Figure \ref{fig: bird}.

\subsection*{Crinkle}

\begin{figure} [p]
	\noindent \begin{centering}
		\includegraphics[width=1\linewidth]{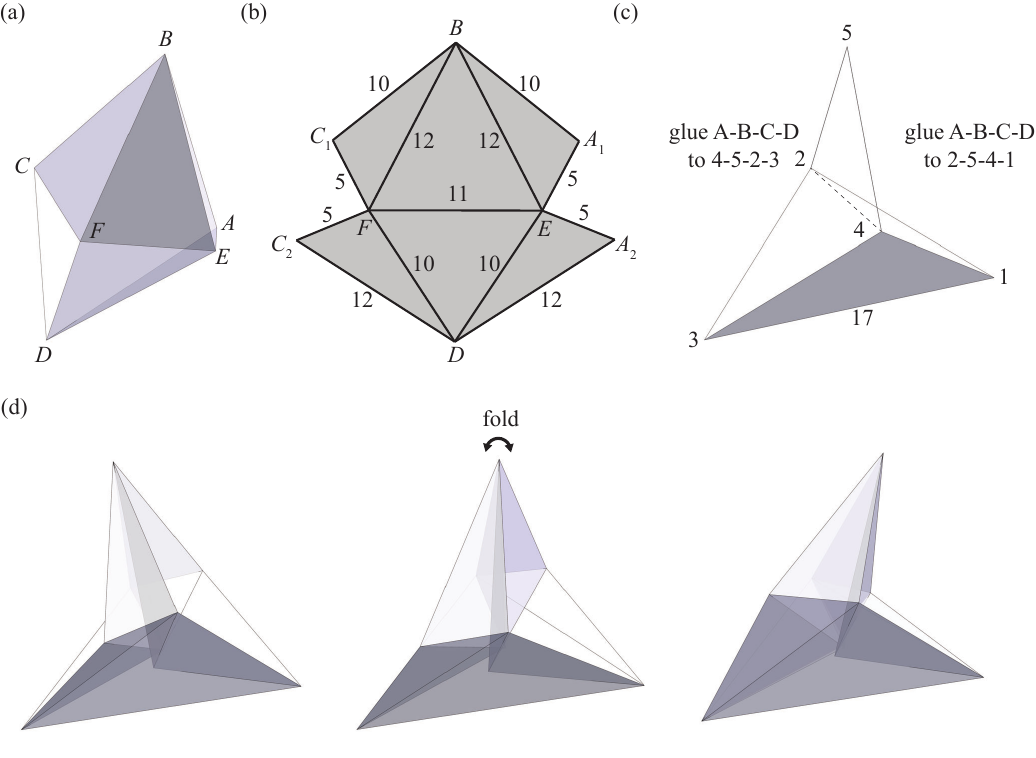}
		\par \end{centering}
	
	\caption{\label{fig: Steffen}(a) A crinkle \(ABCDEF\), obtained by removing two triangular faces from Bricard’s first type of self-intersecting flexible octahedron. (b) A net of (a), with all edge lengths labeled. During the folding process from (b) to (a), the vertex pairs \(A_1, A_2\) and \(C_1, C_2\) are respectively identified with vertices \(A\) and \(C\). Visually, vertex \(E\) `pops in', while vertex \(F\) `pops out' in the resulting three-dimensional configuration (a). (c) A base designed to support the assembly of crinkles, which leads to the construction of Steffen’s polyhedron. Vertices \(1\text{-}2\text{-}3\text{-}4\) define a fixed tetrahedral chamber. The distance between vertices \(1\) and \(3\) is labelled. Triangle \(2\text{-}4\text{-}5\), with fixed edge lengths, can rotate freely about edge \(2\text{-}4\). By attaching two crinkles to the two spatial quadrilaterals \(4\text{-}5\text{-}2\text{-}3\) and \(2\text{-}5\text{-}4\text{-}1\) in opposing configurations (`pop in' vs. `pop out'), we obtain Steffen’s polyhedron (d). The motion of Steffen’s polyhedron closely resembles the motion of the base.}
\end{figure}

Let $ABCD$ be a spatial quadrilateral with fixed edge lengths and a fixed diagonal length $AC$. We define a \textit{crinkle} as a non-self-intersecting flexible polyhedral surface bounded by the quadrilateral $ABCD$, which does not include either of the triangular faces $ABC$ or $ACD$. A crinkle, together with triangles $ABC$ and $ACD$, form a self-intersecting, spherical flexible polyhedron with zero volume. The concept of a crinkle (see Figure~\ref{fig: Steffen}(a)) was introduced by \citet{connelly_counterexample_1977}, who constructed such a surface by removing two triangular faces from each of the three types of Bricard’s self-intersecting flexible octahedra.

A \textit{net} of a polyhedral surface in $\mathbb{R}^3$ is a connected, planar arrangement of polygons obtained by cutting along a subset of the edges of the polyhedral surface so that the result can be unfolded into the plane without overlap. Figure \ref{fig: Steffen}(b) shows a net of \ref{fig: Steffen}(a). 

\subsection*{Base}

A \textit{base} is a structural framework consisting of: (1) a finite set of rigid polyhedral faces and (2) a finite set of bars, each represented as a line segment connecting two vertices. It serves as the core component for assembling non-self-intersecting flexible polyhedra through the attachment of crinkles. In Figure~\ref{fig: Steffen}(c), we illustrate the base used to construct the Steffen's polyhedron and the procedure for assembling it. The kinematic behaviour of the resulting flexible polyhedron closely follows the motion permitted by the base. Since a crinkle contributes no enclosed volume, the volume of the final flexible polyhedron is entirely determined by the base. To ensure non-self-intersection in the assembled structure, the base must be carefully designed to accommodate both geometric compatibility and collision avoidance.

\begin{figure} [tbph]
	\noindent \begin{centering}
		\includegraphics[width=0.9\linewidth]{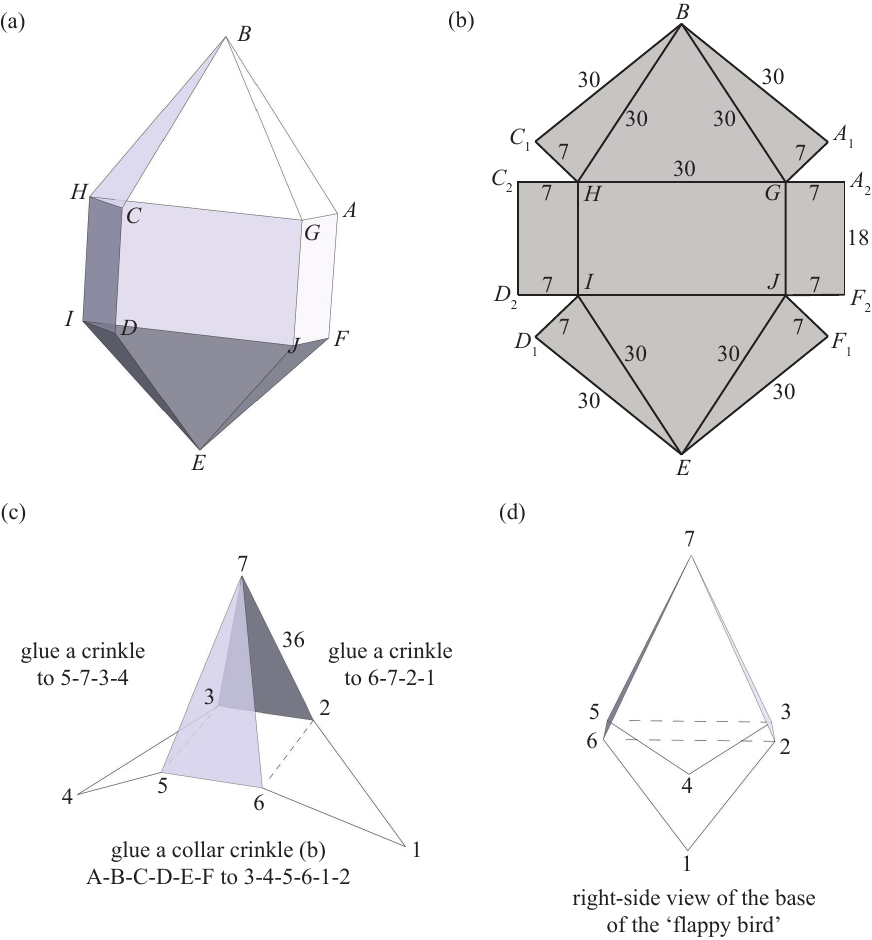}
		\par \end{centering}
	
	\caption{\label{fig: collar crinkle}(a) A collar crinkle \(ABCDEFGHIJ\), obtained by adding a row of rectangles to the crinkle generated from Bricard’s second type of self-intersecting flexible octahedron. (b) A net of (a), with all edge lengths labelled. During the folding process from (b) to (a), the vertex pairs \(A_1, A_2\), \(C_1, C_2\), \(D_1, D_2\), and \(F_1, F_2\) are respectively identified with vertices \(A\), \(C\), \(D\), and \(F\). Visually, vertex \(G\) and \(J\)  `pop in', while vertex \(H\) and \(I\) `pop out' in the resulting three-dimensional configuration (a). (c) A base designed to support the assembly of a collar crinkle and two crinkles, which leads to the construction of Figure \ref{fig: bird}(a). Vertices \(2\text{-}3\text{-}5\text{-}6\text{-}7\) define a fixed right rectangular pyramid. The distance between vertices \(2\) and \(7\) is labelled. Triangles \(1\text{-}2\text{-}6\) and \(3\text{-}4\text{-}5\), with fixed edge lengths, can rotate freely about edge \(2\text{-}6\) and \(3\text{-}5\). By attaching two crinkles to the two spatial quadrilaterals \(6\text{-}7\text{-}2\text{-}1\) and \(5\text{-}7\text{-}3\text{-}4\) in aligned configurations (`pop in' vs. `pop in'), and attaching a collar crinkle to the hexagon \(3\text{-}4\text{-}5\text{-}6\text{-}1\text{-}2\), we obtain the `flappy bird' Figure \ref{fig: bird}(a). The motion of the flappy bird closely resembles the motion of the base. (d) is the right-side view of (c).}
\end{figure}

\subsection*{Collar crinkle}

The idea of `base + crinkle' motivates us to generalize the idea of crinkles. In Figure~\ref{fig: collar crinkle}(a), let \( ABCDEF \) be a spatial hexagon with fixed edge lengths and fixed diagonals \( AC \) and \( DF \), where the quadrilateral \( ACDF \) is required to lie in a plane. We define a \textit{collar crinkle} as a non-self-intersecting flexible polyhedral surface bounded by the hexagon \( ABCDEF \), which does not include either of the faces $ABC$, $ACEF$ or $DEF$. A collar crinkle, together with triangles $ABC$, $DEF$ and quadrilateral $ACDF$, form a self-intersecting, topologically spherical flexible polyhedron with zero volume. Figure \ref{fig: collar crinkle}(b) shows a net of \ref{fig: collar crinkle}(a).

Building on this, we identified a suitable base for the collar crinkle, illustrated in Figure~\ref{fig: collar crinkle}(c). By assembling the collar crinkle at the bottom and attaching two additional crinkles on the sides, we construct a new non-triangulated, non-self-intersecting flexible polyhedron. This new surface exhibits a broader range of motion compared to Steffen’s polyhedron (full polyhedron shown in Figure \ref{fig: bird}(b)), and offers considerable geometric potential as a distortion-free metamorphic gripper.

\subsection*{Assembly into a torus}

Our goal here is to construct a flexible polyhedral torus while minimizing the fixed part. To achieve this, we utilize two nested, topologically spherical flexible polyhedra, ensuring that both the outer surface and the handle remain as flexible as possible. As an example, we designed the base shown in (a) using a nested structure composed of two `flappy bird' bases (Figure~\ref{fig: other examples}(a)). First, a smaller `flappy bird' base is positioned appropriately inside the larger one to avoid self-intersection. Second, by exploiting the common fixed region shared by the two `flappy birds,' we add connecting fixed tunnels to create the handle, thereby assembling the entire structure into a polyhedral torus. The tunnel dimensions are carefully selected to prevent additional self-intersection that may arise during the motion of both nested `flappy birds'. The resulting torus exhibits two kinematic degrees of freedom.

\subsection*{New insight for robotic design}

\begin{figure} [tbph]
	\noindent \begin{centering}
		\includegraphics[width=1\linewidth]{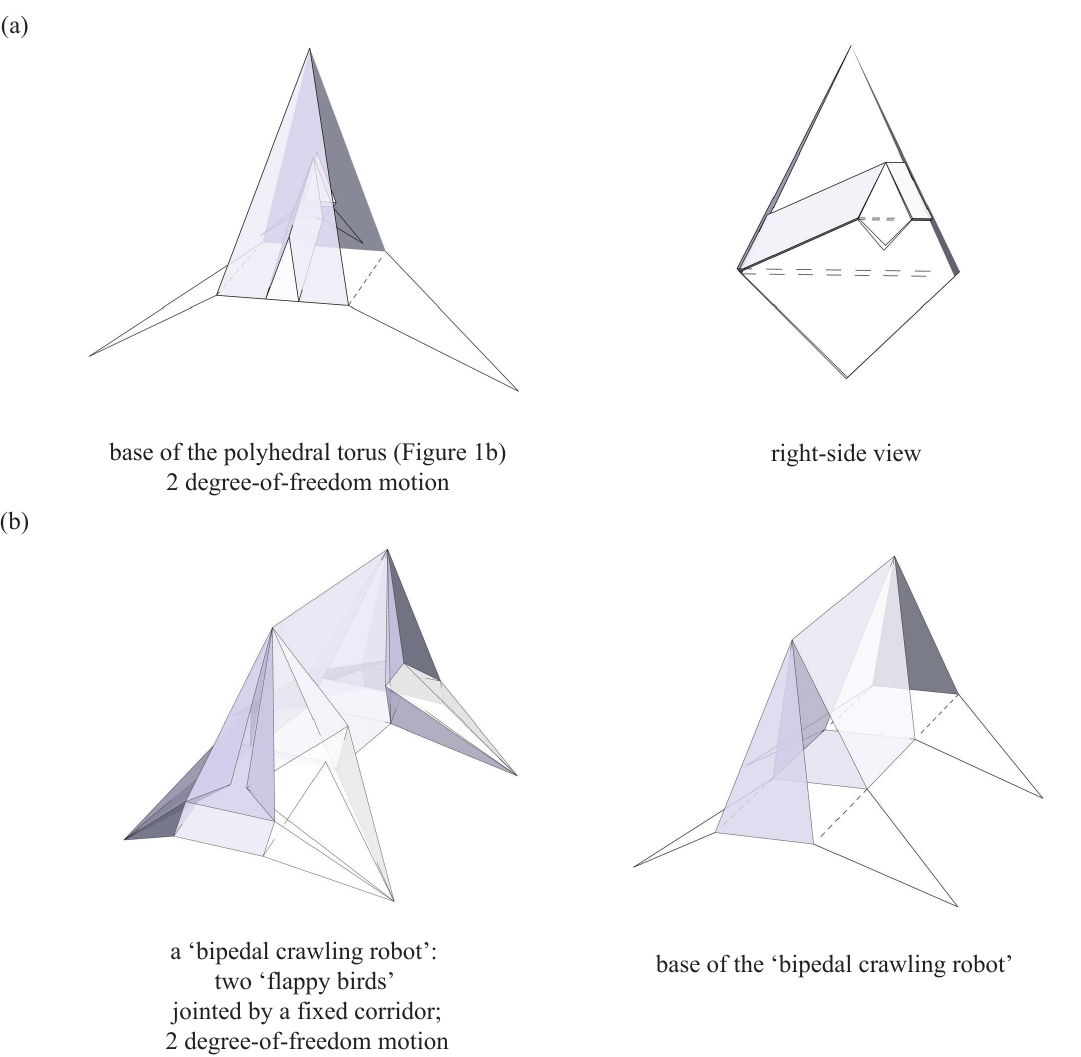}
		\par \end{centering}
	
	\caption{\label{fig: other examples}Further examples generated using the `base + crinkle' method, showing the form-finding potential. (a) shows the base of a polyhedral torus whose full construction is provided in Figure \ref{fig: bird}(b). (b) A topologically spherical flexible polyhedron formed by joining two copies of 'flappy birds' using two fixed trapezoids and one fixed rectangle, with two triangular faces removed. This polyhedron possesses two kinematic degrees of freedom and demonstrates potential applications in crawlable robotics with constant enclosed volume while keeping all faces perfectly rigid.}
\end{figure}

In Figure~\ref{fig: other examples}(b), we present a `bipedal crawling robot' derived from the `flappy bird' structure by reusing its fixed components. The left and right `feet' can move independently and are capable of crawling, for instance, along the ground or up a ladder when equipped with suitable control modules. Throughout its motion, all panels remain perfectly rigid, and the enclosed volume remains constant.

\section{Discussion}

We believe that the `base + crinkle' method provides a promising foundation for the generalization of flexible polyhedra, especially for those with topologies extending beyond the sphere.

\subsection*{Generalization of crinkles}

The current exploration of crinkles and collar crinkles remains far from exhaustive. Even Bricard’s classical results on flexible octahedra are not entirely complete. For crinkles, numerous possibilities extend beyond assemblies of six triangles, yet these have not been systematically investigated. Notably, \citet{atlason_constructing_2025} and \citet{atlason_cutting_2025} proposed two complementary generalizations inspired by the rotational and reflectional symmetries characteristic of Bricard types I and II. The first introduces new triangulated crinkles bounded by quadrilaterals that inherit these respective symmetries. The second constructs larger composite crinkles by cutting matching quadrilateral openings on two flexible polyhedra and joining them along the boundary.

For collar crinkles, we have not yet identified further variants derivable from Bricard I or II under the constraint $BA=BC$, aside from the example illustrated in Figure~\ref{fig: collar crinkle}(a). A broader investigation is needed into how collar crinkles can be systematically derived from various crinkles. Additionally, while our current examples use a hexagonal boundary to facilitate the inclusion of quadrilateral faces, alternative boundary types are entirely feasible and merit further exploration.

Furthermore, there is considerable potential to extend the concept of crinkles beyond constructions derived solely from the symmetries identified in Bricard’s octahedra. Bricard’s octahedra are self-intersecting flexible quadrilateral bipyramids, recent studies have broadened this landscape. \citet{gallet_pentagonal_2024} classified the motions of flexible, self-intersecting pentagonal bipyramids, introducing a new family of promising candidates for crinkle formation. Despite that all flexible $n$-gonal bipyramids must be self-intersected \citep{connelly_attack_1975}, further exploration along this line may yield additional insights—for instance, examining crinkles from flexible $n$-gonal bipyramids that require less unfolded edges when being assembled to suitable bases. Other notable advances include the discovery of a flexible self-intersecting icosahedron \citep{brakhage_icosahedra_2020}, as well as the development of generalized quad-mesh origami mechanisms that resemble collar crinkles \citep{he_rigid_2020}. All of them have different motions compared to Bricard's and provide fertile ground for the development of new classes of crinkles.

\subsection*{High genus flexible polyhedra}

Our current torus construction is achieved by connecting two nested, topologically spherical polyhedra through a compatible cylindrical handle. The handle may either exploit the common fixed region shared by the two polyhedra (Figure \ref{fig: other examples}(a)) or, more complicated, exhibit coordinated motion compatible between them. A promising direction for further exploration is the design of origami tubes that maintain boundary compatibility between specified faces of the two topologically spherical polyhedra. This approach can be naturally generalized to generate flexible polyhedra of higher genus, providing a systematic pathway toward more complex topological configurations.

\subsection*{Motion planning}

Although no unified framework currently exists for designing bases that accommodate diverse input crinkles, once a candidate base is identified, it is natural to consider optimizing its geometric parameters to meet specific performance objectives. Such optimization aims to adjust the base and crinkle dimensions to prevent self-intersection, extend the folding range, and increase the enclosed volume. These objectives can be effectively pursued through numerical techniques, which offer a systematic means of exploring and identifying high-performing configurations for engineering applications.

\subsection*{The function of symmetry}

The role of symmetry in enabling flexibility is fundamental; however, excessive symmetry can lead to self-intersection, and thus a flexible polyhedron cannot be overly symmetric. For example, in Figure \ref{fig: Steffen}(b), if we assign equal edge lengths --- such as setting the lengths of $BE$, $BF$, $DC_2$, $DA_2$ to $10$ --- then vertices $E$ and $F$ would lie symmetrically in the middle plane of $BD$. This symmetry would lead to inevitable self-intersection in the assembled configuration shown in Figure \ref{fig: Steffen}(d). Hence, achieving flexibility while avoiding collisions often requires breaking perfect symmetry and balancing geometric constraints more subtly.

\section*{Acknowledgement}

This work was supported by the CUED–MathWorks Grant. We gratefully acknowledge Elvar Wang Atlason for his insightful contributions during his undergraduate summer project in 2022 and for the stimulating follow-up discussions that ensued. His participation was supported by the Cambridge Mathematics Summer Research Programme, hosted by the Faculty of Mathematics at the University of Cambridge.

We also thank Robert Connelly and Jan Legerský for their inspiring discussions during the Geometry of Materials, Packings, and Rigid Frameworks semester program at the Institute for Computational and Experimental Research in Mathematics (ICERM), Providence, RI, USA. 

\bibliographystyle{plainnat}

\end{document}